\documentclass{article}

\title{A translation of Maehara's ``Eine Darstellung der Intuitionistischen Logik in der Klassischen''\footnote{Translated with permission from Cambridge University Press \\ Originally published in \textit{Nagoya Mathematical Journal}, 7:45-64 in 1954. \\ DOI: \url{https://doi.org/10.1017/S0027763000018055} The present English translation is by Justus Becker. \\ \copyright~1954 Cambridge University Press} }
\author{}
\date{}

\input{macropack}

\begin{document}
\renewcommand{\abstractname}{Abstract by translator}

\maketitle
\begin{abstract}
    A key motivation for Heyting's intuitionistic logic was to gain a formal notion of Brouwer's idea of mathematics as a ``construction of the mind''~{[Hey30]}, one might thus argue that Heyting's Calculus should also correspond to a notion of provability.
    Inspired by this idea, Gödel formalised this connection via an embedding into a modal calculus, which is now known as the modal logic $\Sfour$~{[Göd33]}. 
    While in his original publication, he only proved soundness for the embedding from intuitionistic propositional logic into $\Sfour$, the converse was proved fifteen years later by McKinsey and Tarski~{[MKT48]}. 
    Although, it was later discovered that Gödel also had obtained a proof of the faithfulness of his embedding in unpublished notes in 1941 (see~{[vP25]}).
    Rasiowa and Sikorski later extended Gödel’s embedding to first-order intuitionistic logic~{[RS63]}. 
    In 1954, Maehara independently obtained the same results using proof-theoretic methods, even extending the embedding to one from intuitionistic first-order logic into \emph{intuitionistic} first-order modal logic. 
    This document presents a faithful English translation of Maehara's 1954 paper.
\end{abstract}
\thispagestyle{empty}

\pagebreak
\begin{center}
    \textbf{References for Abstract}
\end{center}
\begin{enumerate}
    \item[{[Göd33]}] Kurt Gödel. ``Eine Interpretation des intuitionistischen Aussagenkalküls''. In: \textit{Ergebnisse eines mathematischen Kolloquiums} 4 (1933), pp. 39–40.
    \item[{[Hey30]}] A.\ Heyting. \textit{Die formalen Regeln der intuitionistischen Logik I}. Sitzungsberichte der Preussischen Akademie der Wissenschaften zu Berlin, 1930, pp. 42–56.
    \item[{[MKT48]}] J.\ C.\ C.\ McKinsey and Alfred Tarski. ``Some Theorems About the Sentential Calculi of Lewis and Heyting''. In: \textit{The Journal of Symbolic Logic} 13.1 (1948), pp. 1–15.
    \item[{[vP25]}] Jan von Plato. “Gödel’s modal interpretation of intuitionistic logic and its proof theory”. In: \textit{Monatshefte für Mathematik} 208.4 (Dec. 2025), pp. 791–817. 
    \item[{[RS63]}] H.\ Rasiowa and R.\ Sikorski. \textit{The Mathematics of Metamathematics}. Monografie matematyczne. PWN-Polish Scientific Publishers, 1963.
\end{enumerate}
\thispagestyle{empty}

\pagebreak

\titleformat{\section}{\bfseries \large}{\S \thesection \ }{0.5ex}{\centering}[]
\setcounter{section}{0}
\titleformat{\subsection}{\textit \normalsize}{\hspace{4mm} \thesubsection. \ }{}{}[]

\begin{center} \setcounter{page}{1}
    \Large \textbf{A Representation of Intuitionistic Logic inside Classical Logic}
\end{center}

\begin{center}
   \large \textsc{Sh\^{o}ji Maehara}
\end{center}
\makeatletter\def\Hy@Warning#1{}\makeatother
{\let\thefootnote\relax\footnotetext{Received on the 20 Oct.\ 1953}}

It might seem inconceivable for anyone who is used to \emph{classical logic} that the ``law of excluded middle'' is generally not recognised in \emph{intuitionistic logic}. However, this seems fairly natural if we observe that a formula, although syntactically the same, may have different contents in classical and intuitionistic logic.

In classical logic, a statement ``$\!A$'' has the intended meaning of ``$\!A$ is in fact true'' -- whether it can be known or not -- however in intuitionistic logic it has the intended meaning of ``it is provable that $A$ is true''. 
The statement ``not $A$'' merely means ``$\!A$ is false'' in classical logic, but ``it is provable that $A$ is false'' in intuitionistic logic.

But, we can also explain the two kinds of logics by providing a translation of every statement in one logic into a statement of similar meaning in the other. Kuroda~\cite{Kuroda_1951}\footnote{See the literature list at the end.} translated every statement of classical logic into a statement of intuitionistic logic. 
In the following, we will make conversely clear that every statement of intuitionistic logic is translatable into a statement of classical logic.

In this treatise, we intend to represent classical and intuitionistic logic by the use of Gentzen's\footnote{See Gentzen \cite{Gentzen1935}} calculi $\LK$ and $\LJ$, as well as related systems; the following investigations therefore relate to the field of \emph{predicate logic}. In the following, we will be using the terminology and symbols after Gentzen~\cite{Gentzen1935}, with few exceptions.

In \ref{sect:preview}, we take a look at the calculus $\LJ$ in preparation of \ref{sect:thm}.

In \ref{sect:Bew}, we introduce a new operator on formulas 
%(Aussagenverknüpfungszeichen) 
``$\Prov$'', i.e.: if $\A$ is a formula, we also view
\[
\Prov \; \A
\]
as a formula, which we read as ``$\A$ is provable''. 

Naturally, $\Prov$ introduced in the calculus $\LK$ means ``provable in classical logic''  -- thus possibly via the law of excluded middle --, and $\Prov$ introduced in the calculus $\LJ$ means ``provable in intuitionistic logic'' -- without the law of excluded middle -- .

From now on, we call $\BLK$ and $\BLJ$ the deductive calculi which are obtained by introducing $\Prov$ to the calculi $\LK$ and $\LJ$ respectively. 
 
In \ref{sect:thm}, we prove the following theorem in which $\A$ is a formula in the calculus $\LJ$, and $\A^b$ is the formula in the calculus $\BLK$ obtained by substituting every subformula $\T$ of $\A$ by $\Prov \T$.\footnote{See \myhyperlink{3.1} for details.}

\textit{Main theorem\footnote{See \myhyperlink{3.7}.}:} $\A$ is derivable in $\LJ$ if and only if $\A^b$ is derivable in $\BLK$.

We can see by this main theorem that a statement, i.e.~a statement that is expressed by a formula $\A$ in intuitionistic logic, and the statement that is expressed by $\A^b$ in classical logic must have the same meaning. 
Nothing else is meant by the title of this treatise. 

Further, we prove the following proposition.

\textit{Corollary\footnote{See \myhyperlink{3.51}.} of theorem \myhyperlink{2}:} Let $\A$ be a formula in the calculus $\LJ$. $\A^b$ is derivable in $\BLJ$ if and only if $\A$ is derivable in $\LJ$.

From this corollary we can see that the concept of \textit{provability} is already included in any proposition in intuitionistic logic.

These two theorems show that for a proof of a formula $\A^b$ we are allowed to use the law of excluded middle, whether one takes it to be valid or not. This is because one could also remove the use of the law of excluded middle from such a proof, if necessary.\footnote{See \myhyperlink{3.6} Theorem \myhyperlink{3}.}

For convenience, we include the \textit{inference figure schemata} of the calculi $\LK$ and $\LJ$. For details, see Gentzen \cite{Gentzen1935}.

Schemata for structural inference figures: 

\renewcommand{\arraystretch}{2.5}
\vspace{2mm}
\begin{tabular}{lc}
    \hspace{-7mm}Thinning: &   $\vlinf{}{}{\D ,\nGamma \to \Theta}{\phantom{D ,} \Gamma \to \Theta} \quad \vlinf{}{}{\Gamma \to \Theta , \D }{\Gamma \to \Theta \phantom{, D}}$ \\
    \hspace{-7mm}Contraction: \hspace{8mm} \phantom{.}& $\vlinf{}{}{\phantom{\D, } \D ,\Gamma \to \Theta}{\D, \D, \Gamma \to \Theta} \quad \vlinf{}{}{\Gamma \to \Theta , \D \phantom{, \D}}{\Gamma \to \Theta , \D, \D}$ \\
    \hspace{-7mm}Exchange: & $\vlinf{}{}{\Delta, \E , \D , \Gamma \to \Theta}{\Delta,  \D , \E ,\Gamma \to \Theta} \quad \vlinf{}{}{\Gamma \to \Theta , \D , \E , \Lambda}{\Gamma \to \Theta, \E, \D, \Lambda}$ \\
    \hspace{-7mm}Cut: & $\vliinf{}{}{\Gamma, \Delta \to \Theta , \Lambda}{\Gamma \to \Theta, \D}{\D , \Delta \to \Lambda}$
\end{tabular}

\vspace{4mm}
Schemata for operational inference figures:

\vspace{2mm}
\begin{tabular}{llll}
    \hspace{-7mm}$\UES$: & $\vliinf{}{}{\Gamma \to \Theta , \A \wedge \B}{\Gamma \to \Theta , \A}{\Gamma \to \Theta , \B}$ &$ \OEA$: & $\vliinf{}{}{\A \vee \B, \Gamma \to \Theta}{\A , \Gamma \to \Theta }{\B, \Gamma \to \Theta}$ \\
    \hspace{-7mm}$\UEA$: & $\vlinf{}{}{\A \wedge \B , \Gamma \to \Theta}{\phantom{\B \wedge\;}\A , \Gamma \to \Theta} \; \; \vlinf{}{}{\A \wedge \B , \Gamma \to \Theta}{\phantom{\A \wedge\;}\B , \Gamma \to \Theta}$ & $\OES$: & $\vlinf{}{}{ \Gamma \to \Theta, \A \vee \B}{\Gamma \to \Theta, \A \phantom{\; \wedge \B}} \;\; \vlinf{}{}{\Gamma \to \Theta, \A \vee \B}{ \Gamma \to \Theta, \B \phantom{\;\; \wedge \A}}$ \\
    \hspace{-7mm}$\AES$: & $\vlinf{}{}{\Gamma \to \Theta ,\forall \rF \F (\rF)}{\Gamma \to \Theta , \F (\aF) \phantom{\forall \rF}}$ & $\EEA$: & $\vlinf{}{}{\exists \rF \F (\rF), \Gamma \to \Theta}{\phantom{\exists \rF} \F (\aF), \Gamma \to \Theta}$ \\
    \hspace{-7mm}$\AEA$: & $\vlinf{}{}{\forall \rF \F (\rF), \Gamma \to \Theta}{\phantom{\forall \rF} \F (\tF), \Gamma \to \Theta}$ & $\EES$: & $\vlinf{}{}{\Gamma \to \Theta ,\exists \rF \F (\rF)}{\Gamma \to \Theta , \F (\tF) \phantom{\exists \rF}}$ \\
    \hspace{-7mm}$\NES$: & $\vlinf{}{}{\phantom{\A, }\Gamma \to \Theta , \neg\A}{\A, \Gamma \to \Theta \phantom{, \neg\A }}$ & $\NEA$: & $\vlinf{}{}{\neg \A , \Gamma \to \Theta \phantom{, \A}}{\phantom{\neg\A, }\Gamma \to \Theta, \A}$ \\
    \hspace{-7mm}$\FES$: & $\vlinf{}{}{\phantom{\A, }\Gamma \to \Theta , \A \vdash \B}{\A , \Gamma \to \Theta , \B \phantom{\; \,\vdash \A}}$ & $\FEA$: & $\vliinf{}{}{\A \vdash \B , \Gamma , \Delta \to \Theta , \Lambda}{ \Gamma \to \Theta , \A}{\B , \Delta \to \Lambda}$
\end{tabular}

\vspace{4mm}
$\A , \B , \D , \E$ range over arbitrary formulas; $\forall \rF \F (\rF)$ and $\exists \rF \F (\rF)$ can be arbitrary formulas of that form, where then $\F (\aF)$ and $\F (\tF)$ are obtained from $\F (\rF)$ by replacing the bound variable $\rF$ by an arbitrary free variable $\aF$ and an arbitrary term $\tF$ respectively. The symbols $\Gamma, \Delta , \Theta$ and $\Lambda$ range over arbitrary, possibly empty sequences of formulas, separated by commata. 

\textit{Variable condition:} The free variable denoted by $\aF$ is not allowed to occur in the lower sequent of this inference figure.

Further, we have the following restriction for $\LJ$-derivations:

``For every derivation sequent no more than one formula is allowed to occur in the succedent.''

\section[A Preview on the Calculus LJ.]{A Preview on the Calculus $\LJ$.} \label{sect:preview}

We introduce a \textit{helping calculus}, which we call $\mLJ$, differing from the calculus $\LK$ as follows.

The inference figures are formed using the inference figure schemata, however with the following restriction: for the schemata of $\FES$, $\NES$ and $\AES$ nothing may be substituted for $\Theta$; thus, this placeholder is empty.

The new calculus $\mLJ$ is equivalent to the calculus $\LJ$.

1.1.\textit{ The notion of equivalence.}\footnote{Gentzen \cite{Gentzen1935}, Section V, \S  1}

\myhypertarget{1.11}. We introduce the following notion of equivalence between sequents: 

Equal sequents are equivalent.

The sequent
\[
\Gamma \to \A_\nu , ... , \A_1
\]
is equivalent to the sequent
\[
\Gamma \to \A_1 \vee ... \vee \A_\nu \qquad (\nu = 1, 2, ...).
\]
(Writing $\A_1 \vee \A_2 \vee ... \vee \A_\nu$ should mean $((\A_1 \vee \A_2) \vee ...) \vee \A_\nu$.)

1.12. We call two derivations equivalent if the end-sequent of one is equivalent to the end-sequent of the other.

1.13. We call two calculi equivalent if every derivation in one calculus can be transformed into an equivalent derivation in the other. 

1.2. \textit{The proof of equivalence between $\LJ$ and $\mLJ$.}

\myhypertarget{1.21}. \textit{Transformation of an $\mLJ$-derivation into an equivalent $\LJ$-derivation.}

\myhypertarget{1.211}. This goes as follows: One starts by replacing every derivation sequent $\Gamma \to \Theta$ with a sequent $\Gamma \to \Theta^\ast$. If $\Theta$ is non-empty, $\Theta^\ast$ refers to the formulas in $\Theta$ in reversed order connected via $\vee$. If $\Theta$ is empty, let $\Theta^\ast$ be empty.

1.212. Now, we already obtained a tree-like system of sequents. Obviously, the end sequent is already equivalent to the end sequent of the $\mLJ$ derivation. All topmost sequents are of the form $\D \to \D$, thus they are initial sequents of an $\LJ$ derivation. 

The inference figures obtained from inference figures of $\mLJ$ are transformed via the following schemata into parts of an $\LJ$-derivation.

1.212.1. The inference figures: the following inferences are immediately subsumed by $\LJ$-inference figures: Thinning in the antecedent, contraction in the antecedent, exchange in the antecedent, as well as $\UEA$, $\OEA$, $\AEA$, $\EEA$, $\FES$, $\NES$ and $\AES$. (In case of a $\AES$ or a $\EEA$ in $\mLJ$, the $\LJ$-variable condition is fulfilled due to the variable condition in $\mLJ$.)

\myhypertarget{1.212.2}. A cut becomes:
\[
\vliinf{}{}{\Gamma, \Delta \to \Lambda^\ast \vee \Theta^\ast}{\Gamma \to \D \vee \Theta^\ast}{\D , \Delta \to \Lambda^\ast}.
\]

(If $\Theta$ and $\Lambda$ are $\A_\mu , ... , \A_1$ and $\B_\nu, ... , \B_1$ respectively, let $\Lambda^\ast \vee \Theta^\ast$ be $\B_1 \vee ... \vee \B_\nu \vee \A_1 \vee ... \vee \A_\mu$. If $\Lambda$ is empty $\Lambda^\ast \vee \Theta^\ast$ means $\Theta^\ast$. If $\Theta$ is empty $\Lambda^\ast \vee \Theta^\ast$ means $\Lambda^\ast$.)

For empty $\Theta$, the inference figure is a cut in $\LJ$.

Let $\Theta=\A_\mu , ... , \A_1$ be non-empty. The inference figure can be derived as follows: 

\renewcommand*{\arraystretch}{1}
\[ \hspace{-7.5em}
\vlderivation{
\vlin{}{\LJ\textnormal{-cut}}{\Gamma , \Delta \to \Lambda^\ast \vee \Theta^\ast}{
\vlid{}{}{\Gamma \to \D \vee \Theta^\ast \phantom{\Gamma \to  \Delta} \overset{
\begin{matrix}
    . & & \\
    & . & \\
    & & .
\end{matrix} \textnormal{\small and so forth \rlap{in the same way}}
}{\D \vee \Theta^\ast , \Delta \to \Lambda^\ast \vee \Theta^\ast}}{
\vliin{}{}{
\D \vee \A_1 \vee \A_2, \Delta \to \Lambda^\ast \vee \A_1 \vee \A_2}
{\vlin{}{}{\D \vee \A_1 , \Delta \to \Lambda^\ast \vee \A_1 \vee \A_2}{
\vliin{}{}{\D \vee \A_1 , \Delta \to \Lambda^\ast \vee \A_1}{
    \vlin{}{}{\D , \Delta \to \Lambda^\ast \vee \A_1}{\vlhy{\D , \Delta \to \Lambda^\ast}}
}{
    \vlin{}{^{\rlap{$\OES$ \small or identity inference figure}}}{\A_1 , \Delta \to \Lambda^\ast \vee \A_1}{\vliq{}{\rlap{\textnormal{\small possibly multiple thinnings and exchanges}}}{\A_1 , \Delta \to \A_1}{\vlhy{\phantom{\Delta, } \A_1 \to \A_1}}}
}
}}
{\vlin{}{}{\A_2 , \Delta \to \Lambda^\ast \vee \A_1 \vee \A_2}{
\vliq{}{\rlap{ $\noteB$ 
\begin{tabular}{l}
    \small{possibly multiple }  \\
    \small{thinnings and exchanges}
\end{tabular}
}}{\A_2 , \Delta \to \A_2}{\vlhy{\phantom{\Delta, } \A_2 \to \A_2}}}
}
}
}
}
\]

1.212.3. Thinning in the succedent becomes:

\[
\vlinf{}{}{\Gamma \to \D \vee \Theta^\ast}{\Gamma \to \Theta^\ast \phantom{\; \,\vee \D}
}
\]

For an empty $\Theta$, the inference figure becomes $\LJ$-thinning in the succedent. Let $\Theta= \A_\nu , ... , \A_1$ be non-empty. From this we make:

\[
\vlderivation{
\vliiq{}{(\myhyperlink{1.212.2}).}{\Gamma \to \D \vee \A_1 \vee ... \vee \A_\nu}{\vlhy{\Gamma \to \A_1 \vee ... \vee \A_\nu}}{
    \vlin{}{}{\A_1 \to \D \vee \A_1}{\vlhy{\A_1 \to \A_1 \phantom{, \vee \D }}}
}
}\footnote{This means the transformation mentioned in \myhyperlink{1.212.2}.}
\]

\myhypertarget{1.212.4}. A contraction in the succedent becomes:

\[
\vlinf{}{}{\Gamma \to \D \vee \Theta^\ast \phantom{\;\, \vee \D}}{\Gamma \to \D \vee \D \vee \Theta^\ast}
\]

From this we make:

\[
\vlderivation{\vliiq{}{(\myhyperlink{1.212.2}).}{\Gamma \to \D \vee \Theta^\ast}{\vlhy{\Gamma \to \D \vee \D \vee \Theta^\ast}}{\vliin{}{}{\D \vee \D \to \D}{\vlhy{\D \to \D}}{\vlhy{\D \to \D}}}
}
\]

\myhypertarget{1.212.5}. An exchange in the succedent becomes:

\[
\vlderivation{\vlin{}{.}{\Gamma \to \Lambda^\ast \vee \E \vee \D \vee \Theta^\ast}{\vlhy{\Gamma \to \Lambda^\ast \vee \D \vee \E \vee \Theta^\ast}}}
\]

From this we make:

\[
\vlderivation{\vliiq{}{(\myhyperlink{1.212.2}).}{\Gamma \to \Lambda^\ast \vee \E \vee \D \vee \Theta^\ast}{\vlhy{\Gamma \to \Lambda^\ast \vee \D \vee \E \vee \Theta^\ast}}{\vlhy{\Lambda^\ast \vee \D \vee \E \to \Lambda^\ast \vee \E \vee \D}}}
\]

Where we put above $\Lambda^\ast \vee \D \vee \E \to \Lambda^\ast \vee \E \vee \D$ the following derivations \myhyperlink{1.212.51} and \myhyperlink{1.212.52} for this sequent.

\myhypertarget{1.212.51}. If $\Lambda$ is empty, the derivation:

\[
\vlderivation{
\vliin{}{.}{\D \vee \E \to \E \vee \D}{\vlin{}{}{\D \to \E \vee \D}{\vlhy{\D \to \D}}}{\vlin{}{}{\E \to \E \vee \D}{\vlhy{\E \to \E}}}
}
\]

\myhypertarget{1.212.52}. If $\Lambda$ is non-empty, the derivation: 

\[
\vlderivation{
\vliin{}{.}{\Lambda^\ast \vee \D \vee \E \to \Lambda^\ast  \vee \E \vee \D}{
    \vliin{}{}{\Lambda^\ast \vee \D \to \Lambda^\ast  \vee \E \vee \D}{
        \vlin{}{}{\Lambda^\ast \to \Lambda^\ast  \vee \E \vee \D}{
        \vlin{}{}{\Lambda^\ast \to \Lambda^\ast  \vee \E}{\vlhy{\Lambda^\ast \to \Lambda^\ast}}
        }
    }{
        \vlin{}{}{\D \to \Lambda^\ast  \vee \E \vee \D}{\vlhy{\D \to \D}
        }
    }
}{
\vlin{}{}{\E \to \Lambda^\ast  \vee \E \vee \D}{
\vlin{}{}{\E \to \Lambda^\ast  \vee \E}{\vlhy{\E \to \E}}
}
}
}
\]

1.212.6. An instance of $\UES$ becomes: 

\[
\vliinf{}{.}{\Gamma \to (\A \wedge \B) \vee \Theta^\ast}{\Gamma \to \A \vee \Theta^\ast}{\Gamma \to  \B \vee \Theta^\ast}
\]

From this we make: 

\[
\vlderivation{
\vliq{}{\noteB
\begin{tabular}{l}
    \small{contractions, exchanges, $\bigl($\myhyperlink{1.212.4}$\bigr)$,}  \\
    \small{and $\bigl($\myhyperlink{1.212.5}$\bigr)$, as far as necessary.}
\end{tabular}}{
\Gamma \to (\A \wedge \B ) \vee \Theta^\ast
}{
\vliiq{}{(\myhyperlink{1.212.2})}{\Gamma , \Gamma \to (\A \wedge \B) \vee \Theta^\ast \vee \Theta^\ast}{\vlhy{\Gamma \to \B \vee \Theta^\ast}}{
\vliq{}{\textnormal{exchanges, as far as necessary}}{\B , \Gamma \to (\A \wedge \B) \vee \Theta^\ast}{
\vliiq{}{(\myhyperlink{1.212.2})}{\Gamma , \B \to (\A \wedge \B) \vee \Theta^\ast}{\vlhy{\Gamma \to \A \vee \Theta^\ast}}{
\vliin{}{}{\A , \B \to \A \wedge \B}{
    \vlin{}{}{\A , \B \to \A}{
        \vlin{}{}{\B, \A \to \A}{\vlhy{\phantom{\B, }\A \to \A}}
    }
}{
    \vlin{}{}{\A , \B \to \B}{\vlhy{\phantom{\A, } \B \to \B}}
}
}
}
}
}
}
\]

1.212.7. An instance of $\OES$ becomes:

\[
\vlinf{}{}{\Gamma \to \A \vee \B \vee \Theta^\ast}{\Gamma \to \A \vee \Theta^\ast \phantom{, \vee \B}}
\]

(the other instance of $\OES$ is done analogously). 

From this we make:

\[
\vlderivation{
\vliiq{}{(\myhyperlink{1.212.2}).}{\Gamma \to \A \vee \B \vee \Theta^\ast}{\vlhy{\Gamma \to \A \vee \Theta^\ast}}{
    \vlin{}{}{\A \to \A \vee \B}{\vlhy{\A \to \A \phantom{,\vee \B} }}
}
}
\]

1.212.8. An instance of $\EES$ becomes:

\[
\vlinf{}{.}{\Gamma \to \exists \rF \F (\rF) \vee \Theta^\ast}{\Gamma \to \F (\tF) \vee \Theta^\ast \phantom{\exists \rF}}
\]

From this we make: 

\[
\vlderivation{\vliiq{}{(\myhyperlink{1.212.2}).}{{\mathit \Gamma} \to \exists \rF \F (\rF) \vee \Theta^\ast}{\vlhy{\Gamma \to \F (\tF) \vee \Theta^\ast }}{
\vlin{}{}{\F (\tF) \to \exists \rF \F(\rF)}{\vlhy{\F (\tF) \to  \F(\tF) \phantom{\exists \rF}}}
}}
\]

1.212.9. An instance of $\NEA$ becomes:

\[
\vlinf{}{.}{\neg \A , \Gamma \to \Theta^\ast \phantom{,\vee \A}}{\phantom{\neg \A , } \Gamma \to \A \vee \Theta^\ast}
\]

From this we make:

\[
\vlderivation{
\vliq{}{\noteB
\begin{tabular}{l}
    \small{exchanges, as}  \\
    \small{far as necessary.}
\end{tabular}}{\neg \A , \Gamma \to \Theta^\ast }{
\vliiq{}{(\myhyperlink{1.212.2})}{ \Gamma, \neg \A \to \Theta^\ast }{\vlhy{ \Gamma \to \A \vee \Theta^\ast }}{
    \vlin{}{}{\A, \neg \A  \to \phantom{\A}}{
    \vlin{}{}{\neg \A , \A \to \phantom{\A}}{\vlhy{\phantom{\neg \A ,} \A \to \A}}
    }
}
}
}
\]

1.212.10. An instance of $\FEA$ becomes:

\[
\vliinf{}{.}{\A \vdash \B , \Gamma , \Delta \to \Lambda^\ast \vee \Theta^\ast}{\Gamma \to \A \vee \Theta^\ast}{\B , \Delta \to \Lambda^\ast}
\]

From this we make:

\[
\vlderivation{
\vliq{}{\noteB
\begin{tabular}{l}
    \small{exchanges, as}  \\
    \small{far as necessary.}
\end{tabular}}{\A \vdash \B , \Gamma , \Delta \to \Lambda^\ast \vee \Theta^\ast}{
\vliiq{}{(\myhyperlink{1.212.2})}{\Gamma ,\A \vdash \B , \Delta \to \Lambda^\ast \vee \Theta^\ast}{\vlhy{\Gamma \to \A \vee \Theta^\ast}}{
    \vlin{}{}{\A, \A \vdash \B , \Delta \to \Lambda^\ast}{
    \vliin{}{}{\A \vdash \B , \A, \Delta \to \Lambda^\ast}{\vlhy{\A  \to \A}}{\vlhy{\B , \Delta \to \Lambda^\ast}}
    }
}
}
}
\]

This completes the transformation of an $\mLJ$-derivation into an equivalent $\LJ$-derivation.

\myhypertarget{1.22}. \textit{Transformation of an $\LJ$-derivation into an equivalent $\mLJ$-derivation.}

An $\LJ$-derivation is a self equivalent $\mLJ$-derivation, therefore no transformation is needed.

Taking the results of \myhyperlink{1.21} and \myhyperlink{1.22} together, it is shown that the calculi $\LJ$ and $\mLJ$ are equivalent.

1.3. \textit{A note on the equivalence of $\mLJ$-sequents.}

As we discussed in \myhyperlink{1.11}, two distinct $\mLJ$ -sequents may still be equivalent. However, in this case they are still \textit{deductively equivalent}.

1.31. If $\Gamma \to \A_1 \vee ... \vee \A_\nu$ is $\mLJ$-derivable, then so is $\Gamma \to \A_\nu , ... ,\A_1$ $\mLJ$-derivable by the following:

\[
\vliinf{}{.}{\Gamma \to \A_\nu , ... ,\A_1}{\Gamma \to \A_1 \vee ... \vee \A_\nu}{\A_1 \vee ... \vee \A_\nu \to \A_\nu , ... ,\A_1}
\]

Where the sequent $\A_1 \vee ... \vee \A_\nu \to \A_\nu , ... ,\A_1$ is $\mLJ$-derivable.

1.32. If $\Gamma \to \A_\nu , ... ,\A_1$ is $\mLJ$-derivable then so is $\Gamma \to \A_1 \vee ... \vee \A_\nu$ (see \myhyperlink{1.21} and \myhyperlink{1.22}).

1.4. In the following investigations, we use $\LJ$ to refer to the calculus $\mLJ$.

\section[New Connective Symbol ``Prov''.]{New Connective Symbol ``$\Prov$''.} \label{sect:Bew}

\myhypertarget{2.1}. We extend the notion of a formula, differing from their definition in the calculi $\LK$ and $\LJ$ as follows:

If $\A$ is a formula then $\Prov \!\A$ is also a formula.

\myhypertarget{2.2}. We add the following schemata for inference figures:

\[
\vlinf{}{\BES}{\Gamma' \to \Prov \A}{\Gamma' \to \A \phantom{\Prov}}
\qquad
\vlinf{}{\BEA}{\Prov \A , \Gamma \to \Theta}{\phantom{\Prov} \A , \Gamma \to \Theta}
.\footnote{Meaning of this designation: \\ $\BES$ (``Provable'' introduction in the succedent) \\ $\BEA$ (``Provable'' introduction in the antecedent).}
\]

In the schema for $\BES$, one can substitute for $\Gamma'$ an arbitrary (possibly empty) sequence of formulas -- separated through commas -- whose outermost connective is $\Prov$.

2.3. We call the calculi that are obtained from $\LK$ and $\LJ$ by adding \myhyperlink{2.1} and \myhyperlink{2.2} $\BLK$ and $\BLJ$ respectively.

2.4. A few remarks on the new schemata for inference figures.

2.41. In the calculi $\BLK$ and $\BLJ$ one can replace the inference figures of $\BEA$ by initial sequents of the following schema:

\[
\Prov \A \to \A.
\]

The equivalence of this initial sequent with the inference figure $\BEA$ is easily shown via cut.

\myhypertarget{2.42}. The inference figures of $\BES$ can be replaced in the calculi of $\BLK$ and $\BLJ$ by initial sequents and inference figures of the following schemata:

\[
\Prov \A \to \Prov( \Prov \A)
\]

and 

\[
\vlinf{}{}{\Prov \Gamma \to \Prov \A}{\Gamma \to \A}
\]

If $\Gamma$ is empty $\Prov \Gamma$ is also empty. If $\Gamma$ is $\A_1 , ... , \A_\nu$ then $\Prov \Gamma$ refers to $\Prov \A_1 , ... , \Prov \A_\nu$ ($\nu= 1, 2, ...$).

The equivalence of these initial sequents and the inference figures with the inference figures of $\BES$ is easily shown via exchanges in the antecedent, cuts and $\BEA$.

2.5. \textit{Examples of derivable sequents in $\BLK$ and $\BLJ$.}

Every sequent of the following form is derivable in both calculi $\BLK$ and $\BLJ$.

2.511. $\Prov (\A \wedge \B) \to \Prov \A \wedge \Prov \B$

\myhypertarget{2.512}. $\Prov \A \wedge \Prov \B \to \Prov (\A \wedge \B)$

2.521. $\Prov ( \Prov \A \vee \Prov \B) \to \Prov \A \vee \Prov \B$

\myhypertarget{2.522}. $\Prov \A \vee \Prov \B \to \Prov ( \Prov \A \vee \Prov \B) $

2.531. $\Prov (\forall \rF \Prov \F (\rF)) \to \Prov (\forall \rF \F(\rF))$

2.532. $\Prov (\forall \rF \F(\rF))\to \Prov (\forall \rF \Prov \F (\rF))$

2.541. $\Prov (\exists \rF \Prov \F (\rF)) \to \exists \rF \Prov \F(\rF)$

\myhypertarget{2.542}. $\exists \rF \Prov \F(\rF) \to \Prov (\exists \rF \Prov \F (\rF))$

\myhypertarget{2.6}. For both the calculus $\BLK$ and the calculus $\BLJ$ the following theorem applies:

Every $\BLK$-derivation can be transformed into a cut-free $\BLK$-derivation with the same final sequent.

This theorem is proved just like Gentzen's proof of his \textit{Hauptsatz}\footnote{Gentzen \cite{Gentzen1935}, III. Section, § 3.} for the calculus $\LK$: via the introduction of ``\textit{mix}'' and two inductions on \textit{grade} and \textit{rank} of a derivation whose last inference figure is a mix and contains otherwise no mix.

For this we have to add the following to obtain the full proof.

2.61. Let rank $= 2$, the mix formula $\M$ occurring in the succedent of the left and in the antecedent of the right upper sequent of the mix as a principal formula where the outermost connective of $\M$ is $\Prov$.\footnote{Gentzen \cite{Gentzen1935}, III. Section, 3.113.3.} The end of the derivation is then:
\[
\vlderivation{\vliin{}{\mix.}{\Gamma' , \Gamma \to \Theta}{
    \vlin{}{\BES}{\Gamma' \to \Prov \A}{\vlhy{\Gamma' \to \A \phantom{\Prov }}}
}{  \vlin{}{\BEA}{\Prov \A , \Gamma \to \Theta}{\vlhy{\phantom{\Prov }\A , \Gamma \to \Theta}}
}}
\]
We transform this to
\[
\vlderivation{\vlin{}{\!\Bigl\{ \!\!\!\!
\begin{tabular}{l}
    \small{possibly multiple thinnings}  \\
    \small{and exchanges.}
\end{tabular}}{\Gamma' , \Gamma \to \Theta}{\vliin{}{\mix}{\Gamma' , \overline{\Gamma} \to \Theta}{\vlhy{\Gamma' \to \A}}{\vlhy{\A , \Gamma \to \Theta}}}}
\]
The sequence of formulas $\overline{\Gamma}$ is obtained from $\Gamma$ by leaving out all occurrences of $\A$.

We can now apply the inductive hypothesis for the grade on the subderivation ending with the sequent $\Gamma' , \overline{\Gamma} \to \Theta$ as it has a smaller grade than the old derivation. Thus, the whole deduction can be transformed into a mix-free deduction. 

2.62. Let rank $> 2$.

\myhypertarget{2.621}. Let the right rank be bigger than 1, the mix formula $\M$ not occurring in the antecedent of the left upper sequent of the mix and the right upper sequent of the mix being a lower sequent of $\BES$.\footnote{Gentzen \cite{Gentzen1935}, III. Section, 3.121.2.} The end of the derivation is then:
\[
\vlderivation{\vliin{}{\mix.}{\Pi, \overline{\Gamma'} \to \overline{\Sigma}, \Prov \A}{\vlhy{\Pi \to \Sigma}}{\vlin{}{\BES}{\Gamma'  \to \Prov \A}{\vlhy{\Gamma'  \to \A \phantom{\Prov}}}}}
\]
The sequences of formulas $\overline{\Sigma}$ and $\overline{\Gamma'}$ are obtained from $\Sigma$ and $\Gamma'$ by removing all occurrences of $\M$.

$\Gamma'$ is a sequence of formulas with $\Prov$ as outermost connective in which the formula $\M$ occurs at least once, thus $\M$ has the form $\Prov \D$.

Let us now introduce two simple helping terms: 

We call two instances of equal formulas \emph{connected} if they occur in the upper sequent(s) and the lower sequent of an inference figure \emph{in accordance} with the inference figure schema -- excluding cut and mix.

The totality of all formula instances in a derivation which we obtain by starting from a single formula instance, taking all of its connected formula instances, then all of the formula instances connected to those, and so on, is called \emph{formula bundle}; we also say: the \emph{bundle} to which the specified formula instance belongs. 

We say: ``A derivation sequent \emph{belongs to the bundle}'' if a formula instance of the sequent belongs to the bundle. 

\phantom{We say:} ``An operational inference figure in the derivation \emph{belongs to the bundle}'', if the principal formula of the inference belongs to the bundle.

\phantom{We say:} ``An inference of thinning \emph{belongs to the bundle}'' if the \emph{thinning formula} belongs to the bundle.

Let us now consider the bundle belonging to the mix formula $\Prov\D (= \M)$ in the succedent of the left upper sequent $\Pi \to \Sigma$.

The sequent which contains at least one topmost formula instance of the bundle is either a lower sequent of $\BES$, a lower sequent of a thinning or an initial sequent.

2.621.1. We begin by replacing every sequent $\Delta \to \Lambda$ belonging to the bundle with a sequent $\Delta , \overline{\Gamma'} \to \Tilde{\Lambda}$. The sequence of formulas $\Tilde{\Lambda}$ is obtained from $\Lambda$ by replacing every formula instance belonging to the bundle by the formula $\Prov \A$.

2.621.2. A $\BES$ belonging to the bundle:
\[
\vlinf{}{}{\Pi_1' \to \Prov \D}{\Pi_1' \to \D \phantom{\Prov}}
\]

becomes:

\[
\vlinf{}{}{\Pi_1' , \overline{\Gamma'} \to \Prov \A}{\Pi_1' \to \D}.
\]

We transform this into
\[
\vlderivation{
\vlin{}{\BES.}{\Pi_1' , \overline{\Gamma'} \to \Prov \A}{
    \vliin{}{\mix}{\Pi_1' , \overline{\Gamma'} \to \Prov \A}{
    \vlin{}{\BES}{\Pi_1' \to \Prov \D}{\vlhy{\Pi_1' \to \D \phantom{\Prov }}}
    }{
    \vlhy{\Gamma' \to \A}
    }
}
}
\]

(We write the corresponding derivation above $\Gamma' \to \A$.)

The derivation of the lower sequent of the new instance of mix has a left rank number of 1, while its right rank number is 1 smaller than the one of the original derivation. 
Therefore, we can eliminate this instance of mix by inductive hypothesis.

2.621.3. A thinning belonging to the bundle:
\[
\vlinf{}{}{\Pi_2 \to \Sigma_2 , \Prov \D}{\Pi_2 \to \Sigma_2 \phantom{, \Prov \D}}
\]

becomes:
\[
\vlinf{}{}{\Pi_2 , \overline{\Gamma'} \to \Tilde{\Sigma}_2, \Prov\A}{\Pi_2 , (\overline{\Gamma'}) \to \Tilde{\Sigma}_2 \phantom{, \Prov \A}}.
\]

($(\overline{\Gamma'})$ means either $\overline{\Gamma'}$ or nothing, depending on whether the sequent $\Pi_2 \to \Sigma_2$ belongs to the bundle or not.)

We replace this by possibly multiple thinnings and exchanges, as far as necessary. 

2.621.4. An initial sequent belonging to the bundle:
\[
\Prov \D \to \Prov \D
\]

becomes:
\[
\Prov \D , \overline{\Gamma'} \to \Prov \A .
\]
From this we make:
\[
\vlderivation{
\vlin{}{\BES.}{\Prov \D , \overline{\Gamma'} \to \Prov \A}{
    \vliin{}{\mix}{\Prov \D , \overline{\Gamma'} \to \A}{
        \vlhy{\Prov \D \to \Prov \D}
    }{
        \vlhy{\Gamma' \to \A}
    }
}
}
\]

(We write the corresponding derivation above $\Gamma' \to \A$.)

Clearly, the condition stated in \myhyperlink{2.2} is fulfilled in $\BES$.

The derivation of the lower sequent of the new instance of mix has a left rank number of 1, while its right rank number is 1 smaller than the one of the original derivation. 
Therefore, we can eliminate this mix by inductive hypothesis.

2.621.5. The end of the derivation becomes:
\[
\vliinf{}{}{\Pi , \overline{\Gamma'} \to \overline{\Sigma}, \Prov \A}{\Pi , \overline{\Gamma'} \to \Tilde{\Sigma}}{\Gamma' \to \Prov \A}
\]

We transform this into
\[
\vliqf{}{\noteB
\begin{tabular}{l}
    \small{possibly multiple con-}  \\
    \small{tractions and exchanges}
\end{tabular}}{\Pi , \overline{\Gamma'} \to \overline{\Sigma}, \Prov \A}{\Pi , \overline{\Gamma'} \to \Tilde{\Sigma}}
\]

2.621.6. All figures that arose from other inference figures are already inference figures of the same kind as the original inference figures. 
Clearly, the variable condition for $\AES$ and $\EEA$ is fulfilled due to the \emph{renaming of free variables}\footnote{Gentzen \cite{Gentzen1935}, III. Section, 3.10.}, whereas the condition stated in \myhyperlink{2.2} is fulfilled for any $\BES$ which does not belong to the bundle.

2.622. Let the right rank number be bigger than 1, let the mix formula not have an occurrence in the antecedent of the left upper sequent of the mix and let the right upper sequent of the mix be a lower sequent of $\BEA$.
This case is essentially done just like all other operational inference figures.

2.623. Let the right rank number be equal to 1. The left rank number is then bigger than 1.

In this case, we cannot have that the left upper sequent of mix is a lower sequent of $\BES$.

Let the left upper sequent of the mix be a lower sequent of a $\BEA$.
This case is essentially done just like all other operational inference figures.

%\pagebreak
\section{Some Theorems} \label{sect:thm}
\myhypertarget{3.1}.
We define in the calculi $\BLK$ and $\BLJ$ a formula $\A^\bw$ which is explained for an arbitrary formula $\A$ inductively as follows.

3.11. If $\A$ is a prime formula:
\[
\A^\bw = \Prov \A.
\]

\vspace{-2em}
3.12. 
\begin{align*}
    (\A \wedge \B)^\bw & = \Prov (\A^\bw \wedge \B^\bw) \\
    (\A \vee \B)^\bw & = \Prov (\A^\bw \vee \B^\bw) \\
    (\A \vdash \B)^\bw & = \Prov (\A^\bw \vdash \B^\bw) \\
    (\neg \A)^\bw & = \Prov (\neg \A^\bw) \\
    (\forall \rF \F (\rF))^\bw & = \Prov (\forall \rF \F^\bw (\rF)) \\
    (\exists \rF \F (\rF))^\bw & = \Prov (\exists \rF \F^\bw (\rF)) \\
    (\Prov \A)^\bw & = \Prov (\Prov \A^\bw) 
\end{align*}
where $\F^\bw (\rF)$ is an abbreviation for $(\F (\rF))^\bw$.

3.2. \emph{Example}:
\[
(A \vee \neg A)^\bw = \Prov (\Prov A \vee \Prov ( \neg \Prov A))
\]
\[
(A \wedge \neg A)^\bw = \Prov (\Prov A \wedge \Prov ( \neg \Prov A))
\]
where $A$ is a propositional variable.

3.3. Every sequent of the following form is derivable in both calculi $\BLK$ and $\BLJ$.

3.31. $\A^\bw \to \Prov \A^\bw$ \hfill (see \myhyperlink{2.42}) \hspace{8em}

3.32. $\A^\bw \wedge \B^\bw \to \Prov (\A^\bw \wedge \B^\bw)$ \hfill (see \myhyperlink{2.512}) \hspace{8em}

3.33. $\A^\bw \vee \B^\bw \to \Prov (\A^\bw \vee \B^\bw)$ \hfill (see \myhyperlink{2.522}) \hspace{8em}

3.34. $\exists \rF \F^\bw (\rF) \to \Prov (\exists \rF \F^\bw (\rF))$ \hfill (see \myhyperlink{2.542}) \hspace{8em}

3.4. Let $\Gamma^\bw$ refer to the sequence of formulas
\[
\A_1^\bw , \A_2^\bw , ... , \A_\nu^\bw
\]
if $\Gamma$ is the sequence of formulas $\A_1 , \A_2 , ... , \A_\nu$ ($\nu =0,1,2, ...$).

\textbf{Theorem \myhypertarget{1}}: If the connective $\Prov$ does not occur in $\Gamma$ and $\Theta$, and if 
\[
\Gamma^\bw \to \Theta^\bw
\]
is $\BLK$- or $\BLJ$-derivable, then
\[
\Gamma \to \Theta
\]
is $\LK$- or, respectively, $\LJ$-derivable.

3.41. \textit{Proof of Theorem \myhyperlink{1}.}

3.411. We firstly define the term ``\textit{removal} of the connective $\Prov$ from an arbitrary formula'' inductively as follows.
The resulting formula from a given formula $\A$ is called $\overline{\A}$.

3.411.1. If $\A$ is a prime formula:
\[
\overline{\A} = \A.
\]

\vspace{-2em}
3.411.2.  
\begin{align*}
    \overline{\A \wedge \B} & = \overline{\A} \wedge \overline{\B} \\
    \overline{\A \vee \B} & = \overline{\A} \vee \overline{\B} \\
    \overline{\A \vdash \B} & = \overline{\A} \vdash \overline{\B} \\
    \overline{\neg \A} & = \neg \overline{\A} \\
    \overline{\forall \rF \F (\rF)} & = \forall \rF \overline{\F (\rF)} \\
    \overline{\exists \rF \F (\rF)} & = \exists \rF \overline{\F (\rF)} \\
    \overline{\Prov \A} & =  \overline{\A} 
\end{align*}

3.412. Change all sequents of the $\BLK$- or, respectively, $\BLJ$-derivation which has $\Gamma^\bw \to \Theta^\bw$ as its end sequent:

From \vspace{-1.3em}
\[
\A_1 , ... , \A_\mu \to \B_\nu , ... , \B_1
\]

make \vspace{-1.3em}
\[
\overline{\A_1} , ... , \overline{\A_\mu} \to \overline{\B_\nu} , ... , \overline{\B_1}.
\]

Now, we already obtained a tree-like system of $\LK$- or, respectively, $\LJ$-sequents. 
The end sequent is obviously $\Gamma \to \Theta$.
The topmost sequents are all of the form $\overline{\D} \to \overline{\D}$, which are all initial sequents of an $\LK$- or, respectively, $\LJ$-derivation. 
The inference figures that were obtained from the $\BLK$- or, respectively, $\BLJ$-inference figures, are all already $\LK$- or, respectively, $\LJ$-inference figures, which are of the same kind as the inference figures from which they were obtained.
However, the possibly existing inference figures that were obtained from $\BES$ or $\BEA$ are identity inference figures.  

Thus, the tree-like system is an $\LK$- or, respectively, an $\LJ$-derivation which has end sequent $\Gamma \to \Theta$.
This proves Theorem \myhyperlink{1}. 

3.5. \textbf{Theorem \myhypertarget{2}}: If
\[
\Gamma \to \Theta
\]
is $\LJ$-derivable then 
\[
\Gamma^\bw \to \Theta^\bw
\]
is $\BLJ$-derivable.

\myhypertarget{3.51}. \textbf{Corollary}: If the connective ``$\Prov$'' does not occur in $\Gamma$ and $\Theta$, then 
\[
\Gamma^\bw \to \Theta^\bw
\]

is $\BLJ$-derivable if and only if
\[
\Gamma \to \Theta
\]
is $\LJ$-derivable.

3.511. Informally, this corollary shows that the connective ``$\Prov$'' is \emph{superfluous} in intuitionistic logic.

3.512. This corollary follows due to Theorems \myhyperlink{1} and \myhyperlink{2}.

3.52. \emph{Proof of Theorem \myhyperlink{2}}.

Take an $\LJ$-derivation which has $\Gamma \to \Theta$ as its end sequent:

3.521. Therefore, there is an $\LJ$-derivation which has $\Gamma \to \Theta^\ast$ (see \myhyperlink{1.211}) as its end sequent such that every sequent occurring in it has at most one formula in its succedent (see \myhyperlink{1.21}).

Call this $\LJ$-derivation $H$. 

3.522. Replace every sequent $\Delta \to \Lambda$ of $H$ by $\Delta^\bw \to \Lambda^\bw$.

Now, we already obtained a tree-like system of $\BLJ$-sequents.

The end sequent is obviously $\Gamma^\bw \to \Theta^{\ast\bw}$. 
The topmost sequents are all of the form $\D^\bw \to \D^\bw$, thus already initial sequents of a $\BLJ$-derivation.
The figures that were obtained by the $\LJ$-inference figures in $H$ get transformed to partial $\BLJ$-derivations.

3.522.1. Every $\LJ$-structural inference figure in $H$ is already transformed into a $\BLJ$-structural inference figure.

3.522.21. A $\AES$ in $H$:
\[
\vlinf{}{}{\Delta \to \forall \rF \F(\rF)}{\Delta \to \F(\aF) \phantom{\forall \rF}}
\]
becomes:
\[
\vlinf{}{}{\Delta^\bw \to \Prov (\forall \rF \F^\bw(\rF))}{\Delta^\bw \to \F^\bw(\aF) \phantom{\Prov(\forall \rF)}}.
\]

From this we make:
\[
\vlderivation{\vlin{}{\BES.}{\Delta^\bw \to \Prov (\forall \rF \F^\bw(\rF))}{
\vlin{}{\AES}{\Delta^\bw \to \forall \rF \F^\bw(\rF)\phantom{\Prov()}}{\vlhy{\Delta^\bw \to \F^\bw(\aF)\phantom{\Prov(\forall \rF)}}}
}}
\]

Clearly, the variable condition is fulfilled in this instance of $\AES$, and the condition mentioned in \myhyperlink{2.2} is fulfilled for $\BES$ (see \myhyperlink{3.1}).

3.522.22. A $\AEA$ in $H$:
\[
\vlinf{}{}{\forall\rF \F(\rF) , \Delta \to \Lambda}{\phantom{\forall \rF} \F(\tF) ,\Delta \to \Lambda }
\]
becomes:
\[
\vlinf{}{}{\Prov(\forall\rF \F^\bw(\rF)) , \Delta^\bw \to \Lambda^\bw}{\phantom{\Prov(\forall\rF)}\F^\bw(\tF) ,\Delta^\bw \to \Lambda^\bw}.
\]

From this we make:
\[
\vlderivation{\vlin{}{\BEA.}{\Prov(\forall\rF \F^\bw(\rF)) , \Delta^\bw \to \Lambda^\bw}{
\vlin{}{\AEA}{\phantom{\Prov()}\forall \rF \F^\bw(\rF) ,\Delta^\bw \to \Lambda^\bw}{\vlhy{\phantom{\Prov(\forall\rF)} \F^\bw(\tF) ,\Delta^\bw \to \Lambda^\bw}}
}}
\]

3.522.23. The transformation works totally analogously for the other connectives $\wedge, \vee, \vdash, \neg$ and $\exists$.

By the aforementioned transformation, we have a $\BLJ$-derivation with end sequent $\Gamma^\bw \to \Theta^{\ast \bw}$, and therefore $\Gamma^\bw \to \Theta^{\ast \bw}$ is $\BLJ$-derivable.

3.523. If $\Gamma^\bw \to \Theta^{\ast \bw}$ is $\BLJ$-derivable then so is $\Gamma^\bw \to \Theta^{\bw}$ $\BLJ$-derivable.

3.523.1. This is clearly the case if $\Theta$ is empty.

3.523.2. If $\Theta$ is non-empty:

To prove this fact, it suffices to show that the sequent:
\[
\Theta^{\ast \bw} \to \Theta^\bw
\]
is $\BLJ$-derivable.

If $\Theta $ is $ \A_\nu , ... , \A_2 , \A_1$ ($\nu = 1,2,...$) then the sequent $\Theta^{\ast \bw} \to \Theta^\bw$ is
\[
\Prov(\Prov(\Prov(\A_1^\bw \vee \A_2^\bw) \vee ... ) \vee \A_\nu^\bw  ) \to \A_\nu^\bw , ... , \A_2^\bw , \A_1^\bw .
\]

This sequent is derivable in the calculus $\BLJ$ via possibly multiple $\OEA$, $\BEA$ and thinnings in the succedent.

This fully proves Theorem \myhyperlink{2}.

\myhypertarget{3.6}. \textbf{Theorem \myhypertarget{3}}: If the sequent
\[
\Gamma^\bw \to \Theta^\bw
\]
is $\BLK$-derivable then it is also $\BLJ$-derivable.

3.61. \textit{Proof of Theorem \myhyperlink{3}.}

If $\Gamma^\bw \to \Theta^\bw$ is $\BLK$-derivable then there is a $\BLK$-derivation which has the sequent $\Gamma^\bw \to \Theta^\bw$ as its end sequent, and in which the inference figure of cut does not occur (see \myhyperlink{2.6}).

We refer to the \emph{grade} of a cut-free $\BLK$ derivation as the number of inference figures $\AES$, $\NES$ and $\FES$ in which the lower sequent has at least two formulas in the succedent.

The grade of a derivation is always either zero or a natural number.

To prove the theorem we do an induction on the grade of the cut-free $\BLK$-derivation.

3.611. In the case of grade = $0$, the $\BLK$-derivation is a $\BLJ$-derivation, and thus: $\Gamma^\bw \to \Theta^\bw$ is $\BLJ$-derivable.

3.612. In the case of grade $> 0$ there is at least one inference figure $\AES, \NES$ or $\FES$ in which the lower sequent has at least two formulas in its succedent.
Let us now consider the formula bundle which belongs to the principal formula $\HF$ of this inference figure (see \myhyperlink{2.621}).

A sequent which contains at least one topmost formula instance of the bundle is either an initial sequent or a lower sequent of an inference figure belonging to the bundle (see \myhyperlink{2.621}). The initial sequent, the thinning and the inference figures for connectives belonging to the bundle have the following form:
\[
\HF \to \HF \qquad \vlinf{}{}{\Gamma_1 \to \Theta_1 , \HF}{\Gamma_1 \to \Theta_1 \phantom{, \HF}} \qquad \vlinf{}{}{\phantom{\Psi,}\, \Gamma_2 \to \Theta_2 , \HF}{\Psi , \Gamma_2 \to \Theta_2 , \Omega}.
\]
Where $\Psi$ and $\Omega$ are either auxiliary formulas of the inference figure or empty, and thus only dependent on the principal formula $\HF$.
(At least one operational inference figure belonging to the bundle has at least two formulas in the succedent of its lower sequent.)

The final sequent of the derivation is $\Gamma^\bw \to \Theta^\bw$ while the derivation does not contain cuts, thus the sequent, which contains the lowest formula instance of the bundle, is an upper sequent of $\BES$ --- we call it the \textit{final inference figure} of the bundle ---.

3.612.1. We firstly replace every sequent $\Delta \to \Lambda$ which belongs to the bundle by a sequent $\Delta , \Psi \to \Omega , \Tilde{\Lambda}$. Whereas the sequence of formulas $\Tilde{\Lambda}$ is obtained from $\Lambda$ by omitting all formulas that belong to the bundle.

3.612.2. An initial sequent belonging to the bundle becomes:
\[
\HF , \Psi \to \Omega .
\]
The outermost connective of $\HF$ is $\forall$, $\neg$ or $\vdash$, thus the new sequent has the form:
\[
\forall \rF \F (\rF) \to \F (\aF) \qquad \neg \A , \A \to  \quad \textnormal{ or } \quad \A \vdash \B , \A \to \B
\]
Then, this sequent is $\BLJ$-derivable.

We write the respective $\BLJ$-derivation above $\HF , \Psi \to \Omega$.

3.612.3. A thinning belonging to the bundle becomes: 
\[
\vlinf{}{}{\Gamma_1 , \Psi \to \Omega , \Tilde{\Theta}_1}{\Gamma_1 , (\Psi) \to (\Omega) , \Tilde{\Theta}_1}.
\]
($(\Psi)$ and $(\Omega)$ either mean $\Psi$ and $\Omega$ or nothing, depending on whether the sequent $\Gamma_1 \to \Theta_1$ belongs to the bundle.)

We replace this by possibly multiple thinnings and exchanges, as far as necessary. 

3.612.4. An operational inference figure belonging to the bundle becomes:
\[
\vlinf{}{}{\Gamma_2 , \Psi \to \Omega , \Tilde{\Theta}_2}{\Psi , \Gamma_2 , (\Psi) \to (\Omega) , \Tilde{\Theta}_2 , \Omega}.
\]
We replace this by possibly multiple contractions and exchanges as far as necessary.   

3.612.5. The final inference figure
\[
\vlinf{}{}{\Gamma_3' \to \Prov \HF}{\Gamma_3' \to \HF \phantom{\Prov}}
\]
becomes:
\[
\vlinf{}{}{\phantom{\Psi, } \Gamma_3' \to \Prov \HF}{\Gamma_3' , \Psi \to \Omega \phantom{\Prov}}.
\]

We make from this:
\[
\vlderivation{
\vlin{}{\BES .}{\phantom{\Psi, }\Gamma_3' \to \Prov \HF}{
\vlin{}{\AES, \NES \text{ or } \FES}{\phantom{\Psi, }\Gamma_3' \to \HF \phantom{\Prov}}{
\vliq{}{\noteB
\begin{tabular}{l}
    \small{exchanges as}  \\
    \small{far as necessary}
\end{tabular}}{\Psi , \Gamma_3' \to \Omega \phantom{\Prov}}{\vlhy{ \Gamma_3', \Psi \to \Omega\phantom{\Prov}}}
}
}
}
\]

In the case of $\AES$ the variable condition is clearly fulfilled due to the renaming of free variables. 

3.612.6. The figures obtained by other inference figures which contain formulas of the bundle are already inference figures of the same kind as the old inference figures, while the figures obtained from exchanges which have formulas of the bundle as exchange formulas are the identical inference figure. 
(We clearly have that a $\BES$ which contains formulas of the bundle is the final inference figure of the bundle, while the inference figures of $\AES$ and $\EEA$, obtained from $\AES$ and $\EEA$ respectively, fulfil the variable condition due to the renaming of free variables.)

By this transformation, we have a cut-free $\BLK$-derivation, which has $\Gamma^\bw \to \Theta^\bw$ as its end sequent, and its grade is smaller than the grade of the original $\BLK$-derivation.
Thus, via the inductive hypothesis $\Gamma^\bw \to \Theta^\bw$ is $\BLJ$-derivable.

This proves Theorem \myhyperlink{3}.

\myhypertarget{3.7}. From the proved Theorems \myhyperlink{1}, \myhyperlink{2} and \myhyperlink{3} it follows

\textbf{Main Theorem}: If the connective $\Prov$ does not occur in $\Gamma$ and $\Theta$, then
\[
\Gamma \to \Theta
\]
is $\LJ$-derivable if and only if
\[
\Gamma^\bw \to \Theta^\bw
\]
is $\BLK$-derivable.

\printbibliography

\textit{University of Tokyo}
\end{document}